 \documentclass[letterpaper, 10 pt, conference]{ieeeconf}  
\newcommand{\fref}[1]{(\ref{eq:#1})}

\IEEEoverridecommandlockouts                              

\usepackage{graphicx}
\title{\LARGE \bf
Model Predictive Load Scheduling Using Solar Power Forecasting
}

\author{Abdulelah H. Habib, Jan Kleissl and Raymond A. de Callafon
\thanks{Abdulelah H. Habib, Jan Kleissl and Raymond A. de Callafon are with the Department of Mechanical and Aerospace Engineering, University of California, San Diego
        {\tt\small Ahhabib,jkleissl,callafon@ucsd.edu}}%
}

\begin{document}

\maketitle
\thispagestyle{empty}
\pagestyle{empty}

\begin{abstract}
In this paper a model is developed to solve the on/off scheduling of (non-linear) dynamic electric loads based on predictions of the power delivery of a (standalone) solar power source. Knowledge of variations in the solar power output is used to optimally select the timing and the combinations of a set of given electric loads, where each load has a desired dynamic power profile. The optimization exploits the desired power profiles of the electric loads in terms of dynamic power ramp up/down and minimum time on/off of each load to track a finite number of load switching combinations over a moving finite prediction horizon. Subsequently, a user-specified optimization function with possible power constraints is evaluated over the finite number of combinations to allow for real-time computation of the optimal timing and switching of loads. A case study for scheduling electric on/off loads with switching dynamics and solar forecast data at UC San Diego is carried out.
\end{abstract}

\section{INTRODUCTION}
Power variability and intermittency are the main obstacles facing renewable energy integration into smart grids to create a sustainable electric power system. This problem is more critical in standalone or island mode applications where a high penetration of renewable power sources may create power variability that is large enough to frequently influence electric power quality and increase operating costs.  Storage systems in different formats of electrochemical, mechanical or thermal storage have been applied to solve this problem \cite{c1}, but this will add cost and complexity to the standalone system. An alternative is load scheduling, where loads on the (standalone) system are scheduled to absorb solar power variability, thereby reducing the need for solar curtailment and associated energy losses. 

Load scheduling  has been applied to many load types, such as thermal loads, residential appliances, power industry, and EV charging \cite{c8}. In \cite{c3}, a case study was implemented to overcome the wind power variability through EV charging. An example for residential appliance scheduling was shown in \cite{c4}. Also HVAC systems have been included in the scheduling problem \cite{c5}. From the supply side for a hybrid system was discussed in \cite{c9}. 
In \cite{c2} game theory and customer effects on the grid and EVs were investigated. 

Most approaches to optimal load or demand scheduling use a form of Model Predictive Control (MPC) \cite{c12,c13} to compute optimal control or scheduling signals for the load. Typically, in MPC a constrained (quadratic) optimization problem is solved iteratively over a finite horizon $N$ and a moving time horizon $t$ from $t=k$ till $t=k+N-1$to compute an optimal control signal in real time, denoted here by $w(k)$ at time $t=k$. Countless examples of innovative MPC based approaches for load scheduling, grid tied storage systems or maintaining voltage stability can be found in e.g., \cite{c3}, \cite{c14}, \cite{c15} or \cite{c18}. Although MPC approaches are extremely powerful in computing optimal control signals over a moving but finite time horizon, typically the control signal $w(k)$ is allowed to attain any real value during the optimization, see e.g., \cite{c16,c17,c19}. Unfortunately, a real-valued control signal $w(k)$ would require electrical loads to operate at fractional load demands. Although fractional loads can be accommodated by electrical energy storage systems or partial or pulse width modulation of loads \cite{c15}, (non-linear) dynamic power profiles of the electric loads during power ramp up/down and minimum time on/off of each load is harder to implement in a standard MPC framework.

In this paper we define load scheduling as the optimal on/off combinations and timing of a set of distinct electric loads via the computation of an optimal binary control signal $w(k) \in \{0,1\}$. The work is partially motivated by previous work \cite{c10} and \cite{c11} in which the design sizing problem of a standalone photovoltaic reverse osmosis (RO) system is considered, where the RO loads are to be scheduled on/off. The work in \cite{c10} computes the optimal size and number of units for a selected location but does not consider optimal scheduling of dynamic loads. Here we aim for finding the optimal load schedule and on/off switching events for possibly non-linear dynamics of electric loads. For the MPC solution over a moving prediction horizon, solar forecasting data is used as an input to our model. 

The solar forecasting model developed by UC San Diego utilizes a ground based sky imager to detect clouds and forecast their movement up to 15~min into the future. Using the forecasts of projected cloud shadows on the surface coupled with a clear sky irradiance model yields the Global Horizontal Irradiance (GHI) forecast for a selected location\cite{c1,c2}. With the finite prediction horizon in MPC it is crucial to have reliable and accurate forecasts of power delivery and it is assumed to be provided by the work in \cite{c6,c7} as a starting point for the dynamic load scheduling in this paper.

For binary load switching, MPC optimization problems often become untractable due to a combinatorial problem where the number of binary combinations grows exponentially with the length of the prediction horizon $N$ and the number $n$ of loads. In this paper it is shown that constraints on the allowable load switching help to alleviate the combinatorial problem, making an MPC optimization with binary switching computationally feasible.
This paper is organized as follows. In Section~\ref{sec:SDLM} the assumptions for the dynamic loads such as dynamic power ramp up/down and minimum on/off time are specified. For brevity, linear dynamics for the load models are used, but the rate of change of the load models are chosen to differ for the "on" or "off" switching of the load. In section~\ref{sec:DLS} the approach for dynamic load scheduling is presented. Based on power tracking over a moving prediction horizon of $N$ points, the admissible set of binary switching combinations is derived. It is shown that the number of binary combinations growths much less than for a typical exponential growth with the length of the prediction horizon and the number of loads due to the imposed constraints on the allowable load switching. In Section~\ref{sec:AEx} an illustrative example of a standalone solar system connected to 3 distinctive loads is presented. Solar forecasting data of both a clear and a cloudy day obtained from the Solar Resource Assessment \& Forecasting Laboratory (SRAF) at UC San Diego will be used to illustrate how loads are scheduled to turn on/off dynamically to track solar power predictions. The paper ends by concluding remarks in Section~\ref{sec:CCLS}.
\section{Switched Dynamic Load Modeling}\label{sec:SDLM}
\subsection{Assumptions on Loads}
We consider a fixed number of $n$ loads where the power demand $p_i(t)$, $i=1,2,\ldots,n$ as a function of time $t$ for each load $i$ is modeled by a known switched dynamic system. For the dynamic scheduling of the loads, we will assume that loads can be switched ``on'' or ``off'' by a binary switching signal $w_i(t) = \{0,1\}$. To allow for a realistic dynamic load switching, the dynamic system used to model the power demand $p_i(t)$ may be {\em distinct\/} for each load $i=1,2,\ldots,n$. Furthermore, the distinct dynamic system for a particular load $i$ may also switch its dynamics, depending on the transition of the binary switching signal $w_i(t) = \{0,1\}$. The switched dynamic system assumption allows different dynamics for the time dependent power demands $p_i(t)$ when the binary switching signal $w_i(t) = \{0,1\}$ transitions from 0 to 1 or  from 1 to 0. 

Each load $i$ is also assumed to have a known minimum duration $T_i^{\mbox{\tiny off}}>0$ for the ``off'' time of the load when $w_i(t)=0$ and a minimum duration $T_i^{\mbox{\tiny on}}>0$ for the ``on'' time of the load when $w_i(t)=1$. The duration times $T_i^{\mbox{\tiny off}}$ and $T_i^{\mbox{\tiny on}}$ avoid unrealistic on/off chattering of the switch signal $w_i(t)$ during load scheduling and limit the number of transitions in $w_i(t)$ over a finite optimization period $T>0$. Finally, it is assumed that all loads are switched ``off'', e.g., $w_i(t)=0$ for $i=1,2,\ldots,n$, outside the finite time interval $t\in[0,T]$. 

\subsection{Admissible Switching Signals}

With the minimum on/off duration times $T_i^{\mbox{\tiny on}},T_i^{\mbox{\tiny off}}$ and the finite time period $T$ for load switching, on/off switching of a load at time $t=\tau_i$ can now be formalized. Special care should be given to turning on loads at $t= \tau_i$ close to the final time $\tau_i=T$. For the formalization, the load switching signal $w_i(t)$ will be a combination of an ``on'' signal $w_i^{\mbox{\tiny on}}(t) \in {0,1}$ and an ``off'' signal $w_i^{\mbox{\tiny on}}(t) \in {0,1} $ that both take into account the constraints of minimum on/off duration and the finite time $T$ for load switching. 

As a result, the admissible on/off transition signal $w_i(t) = \{w_i^{\mbox{\tiny on}}(t),w_i^{\mbox{\tiny off}}(t)\}$ of a load at time $t=\tau_i$ can now be formalized by the switching signal
\begin{equation}
w_i^{\mbox{\tiny on}}(t) = 
\left \{ 
\begin{array}{rclcl} 
0 &\mbox{for}& t < \tau_i & \mbox{and}& \tau_i \geq T_{i,last}^{\mbox{\tiny off}} + T_i^{\mbox{\tiny off}} \\ 
1 &\mbox{for}&  t \geq \tau_i & \mbox{and} & \tau_i \leq T - T_i^{\mbox{\tiny on}}
\end{array} \right .
\label{eq:wion}
\end{equation}
where $T_{i,last}^{\mbox{\tiny off}}$ denotes the most recent (last) time stamp at which the load $i$ was switched ``off'',
and
\begin{equation}
w_i^{\mbox{\tiny off}}(t) = 
\left \{ 
\begin{array}{rclcl} 
1 &\mbox{for}& t < \tau_i & \mbox{and}& \tau_i \geq T_{i,last}^{\mbox{\tiny on}} + T_i^{\mbox{\tiny on}} \\ 
0 &\mbox{for}&  t \geq \tau_i
\end{array} \right .
\label{eq:wioff}
\end{equation}
where $T_{i,last}^{\mbox{\tiny on}}$ denotes the most recent (last) time stamp at which the load $i$ was switched ``on''.

\subsection{Dynamic Load Models}

It is clear that the switching time(s) $\tau_i$ for the signal $w_i(t) = \{0,1\}$ depends on the dynamics of the power demand $p_i(t)$, which again may be different for each load. For the computational results presented in this paper, linear continuous-time dynamic models will be used to model the switching dynamics of the power demand of the loads. It should be pointed out that the computational analysis is not limited to the use of linear dynamic models, as long as the dynamic models allow the numerical computation of power demand $p_i(t)$ as a function of the switching signal $w_i(t)$.

To allow different dynamics for the time dependent power demands $p_i(t)$ when the binary switching signal $w_i(t) = \{0,1\}$ transitions from 0 to 1 ("on") or from 1 to 0 ("off"), different dynamics is sued for the load models. This allows power demands $p_i(t)$ to be modeled at different rates when switching loads. Using the Laplace transform ${\cal L}\{ \cdot \}$ and referring back to the admissible on/off transition signals $w_i^{\mbox{\tiny on}}(t)$ and $w_i^{\mbox{\tiny off}}(t)$ respectively in \fref{wion} and \fref{wioff}, the switched linear order continuous-time dynamic models for the loads are assumed to be of the form
\begin{equation}
p_i(s) = G_{i}^{on}(s) x_i w_i(s)~\mbox{and}~ w_i(s) = {\cal L} \{ w_i^{\mbox{\tiny on}}(t)\}
\label{eq:odeon}
\end{equation}
and
\begin{equation}
p_i(s) = G_{i}^{off}(s) x_i w_i(s)~\mbox{and}~ w_i(s) = {\cal L} \{ w_i^{\mbox{\tiny off}}(t)\}
\label{eq:odeoff}
\end{equation}
where $G_{i}^{on}(s)$ and $G_{i}^{off}(s)$ represent the dynamics of the power flow for turning the load $i$ "on" or "off". Both models satisfy $G_{i}^{on}(0)=1$ and $G_{i}^{off}(0)=1$ and a steady-state load demand parameter $x_i$ is used to model the size of the load, but different dynamics is used to model respectively the on/off dynamic switching of the load. As an example, if the $G_{i}^{on}(s)$ is given by a first order system with only a time constant $\alpha_i^{\mbox{\tiny on}}$, the power demand $p_i(t)$ progresses from the ``off'' state $p_i^{\mbox{\tiny off}}$ (not necessarily zero) to the ``on'' state over time for the ``on'' switching of the load. Based on the solution 
\[
p_i(t) = p_i^{\mbox{\tiny off}} \cdot e^{-t/\alpha_i^{\mbox{\tiny on}}} +  x_i (1-e^{-t/\alpha_i^{\mbox{\tiny on}}})
\]
it is clear that 
\[
\lim_{t \rightarrow \infty} p_i(t) = x_i
\]
Clearly, more complex models $G_{i}^{on}(s)$ and $G_{i}^{off}(s)$ may be used, but in this paper we will restrict ourselves to first and second order models for brevity without loss of generality of the approach presented here.

\subsection{Discretization of Models}

A tractable numerical implementation of the computation of the optimal switching times $\tau_i$ of the binary switching signals $w_i(t)$ for each load can be achieved by discretizing the power demand $p_i(t)$ and the optimal switching signal $w_i(t)$ at a time step
\[
t_k = k \Delta_t
\]
where $\Delta_t$ is the sampling time $k=0,1,\ldots$ is an integer index. To simplify the integer math, we assume that 
both the switching times
\begin{equation}
\tau_i = N_i \Delta_t
\label{eq:Ni}
\end{equation}
and the minimum on/off duration times
\begin{equation}
\begin{array}{rcl}
 T_i^{\mbox{\tiny on}} &=& N_i^{\mbox{\tiny on}} \Delta_t \\
 T_i^{\mbox{\tiny off}} &=& N_i^{\mbox{\tiny off}} \Delta_t  
\end{array}
\label{eq:Nionoff}
\end{equation}
are all multiple of the sampling time $\Delta_t$.

With the imposed time discretization, the switching signal $w_i(t_k)$ is held constant between subsequent time samples and $t_k$ and $t_{k+1}$. In that case, the computation of $p_i(t_k)$ can be achieved using a Zero Order Hold (ZOH) discrete-time equivalent of the continuous-time models given earlier in \fref{odeon} and \fref{odeoff}. Using the z-transform ${\cal Z}\{ \cdot \}$, the ZOH discrete-time equivalent dynamic models are given by
\[
p_i(z) = G_{i}^{on}(z) x_i w_i(z)~\mbox{and}~ w_i(z) = {\cal Z} \{ w_i^{\mbox{\tiny on}}(t_k)\}
\]
for ``on'' switching of the load and 
\[
p_i(z) = G_{i}^{off}(z) x_i w_i(z)~\mbox{and}~ w_i(z) = {\cal Z} \{ w_i^{\mbox{\tiny off}}(t_k)\}
\]

for ``off'' switching of the load, where $G_{i}^{on}(z)$ and $G_{i}^{off}(z)$ are the ZOH discrete-time equivalents of $G_{i}^{on}(s)$ and $G_{i}^{off}(s)$ using a sampling time $\Delta_t$. In this way the (discrete-time) load dynamics is fully determined by $G_{i}^{on}(z)$, $G_{i}^{off}(z)$, static load demand $x_i$ and the chosen sampling time $\Delta_t$. 

\section{Dynamic Load Scheduling}\label{sec:DLS}

\subsection{Power Tracking}

To formulate the optimization that allows for a computation of optimal discrete-time switching signals $w_i(t_k),~ i=1,2,\ldots,n$ for the power demand $p_i(t_k)$ of $n$ loads, we consider the problem of tracking a desired discrete-time power profile $P(t_k)$. Defining a power tracking error
\begin{equation}
e(t_k) = P(t_k) -\sum_{i=1}^n p_i(t_k)
\label{eq:powererror}
\end{equation}
it is clear that computing optimal $w_i(t_k)$ involves a criterion function and possible constraints on $e(t_k)$ and $w_i(t_k)$ over a (finite) time horizon $k=1,2,\ldots,N$. Choosing $N$ to be large, e.g., $N = T/\Delta_t$ where $T$ is the complete optimization period $T>0$, suffers from two major drawbacks. 

The first drawback is that the number of possible combinations of the discretized binary switching signal $w_i(t_k)$ grows exponentially with the number of loads $n$ and the number of time steps $N$. At each time step $t_k$ with $n$ loads, one would typically have $2^n$ binary (on/off) load combinations. Starting with an initial combination at time step $k=1$, evaluating possible switching along the remaining $N-1$ time steps leads to $(2^n)^{N-1}$ binary combinations, which quickly becomes intractable. The second drawback of choosing $N$ to be large requires the discrete-time power profile $P(t_k)$ to be available over many time samples to plan for optimal load scheduling. In case of solar power tracking, this would require predictions of the solar power output $P(t_k)$ many time steps ahead, which may quickly becomes unreliable. If indeed $P(t_k)$ is produced by solar power forecasting, it makes sense to limit the prediction horizon $N$ used in load scheduling.

Fortunately, the effects of the first drawback is significantly reduced by the requirement of minimum on/off duration times $T_i^{\mbox{\tiny on}},T_i^{\mbox{\tiny off}}$ for the loads. As mentioned before, this avoids unrealistic on/off chattering of the switch signal $w_i(t)$ during load scheduling and significantly reduces the number of binary load combinations.

\subsection{Admissible Discrete-Time Switching Combinations}

With the imposed time discretization given in \fref{Ni}, \fref{Nionoff} and a finite prediction horizon $N$, the admissible on/off transition signal in \fref{wion} reduces to
\begin{equation}
w_i^{\mbox{\tiny on}}(t_k) = 
\left \{ 
\begin{array}{rclcl} 
0 &\mbox{for}& k < N_i & \mbox{and}& N_i \geq N_{i,last}^{\mbox{\tiny off}} + N_i^{\mbox{\tiny off}} \\ 
1 &\mbox{for}&  k \geq N_i & \mbox{and} & N_i \leq N - N_i^{\mbox{\tiny on}}
\end{array} \right .
\label{eq:wiondiscrete}
\end{equation}
where $N_{i,last}^{\mbox{\tiny off}}$ now denotes the most recent discrete-time index at which the load $i$ was switched ``off''. Similarly, \fref{wioff} reduces to
\begin{equation}
w_i^{\mbox{\tiny off}}(t_k) = 
\left \{ 
\begin{array}{rclcl} 
1 &\mbox{for}& k < N_i & \mbox{and}& N_i \geq N_{i,last}^{\mbox{\tiny on}} + N_i^{\mbox{\tiny on}} \\ 
0 &\mbox{for}&  k \geq N_i
\end{array} \right .
\label{eq:wioffdiscrete}
\end{equation}
where $N_{i,last}^{\mbox{\tiny on}}$ denotes the most recent discrete-time index at which the load $i$ was switched ``on''. Collectively, the signals $w_i^{\mbox{\tiny on}}(t_k)$ in \fref{wiondiscrete} and $w_i^{\mbox{\tiny off}}(t_k)$ \fref{wioffdiscrete} define a set ${\cal W}$ of binary values for admissible discrete-time switching signals defined by
\begin{equation}
{\cal W} = \left \{ 
\begin{array}{c}
w_i(t_k) \in \{w_i^{\mbox{\tiny on}}(t_k),w_i^{\mbox{\tiny off}}(t_k)\} ,\\
 i=1,2,\ldots,n,~ k=1,2,\ldots,N\\
\mbox{where}~
\begin{array}{c}
w_i^{\mbox{\tiny on}}(t_k) \in \{0,1\}~\mbox{given in}~\mbox{\fref{wiondiscrete}}\\
w_i^{\mbox{\tiny off}}(t_k) \in \{0,1\}~\mbox{given in}~\mbox{\fref{wioffdiscrete}}
\end{array}
\end{array}
\right \}
\label{eq:Wset}
\end{equation}

It is worthwhile to note that the number of binary elements in the set ${\cal W}$ is always (much) smaller than $(2^n)^{N-1}$ due to required minimum number of on/off samples $N_i^{\mbox{\tiny on}},N_i^{\mbox{\tiny off}}$ for the loads. This can be seen by using the simplified binomial formula
\[
\begin{array}{rcl}
2^n &=& \displaystyle \sum_{m=0}^{n} \left ( \begin{array}{c} n \\ m \end{array} \right )\\
 &=& \displaystyle 1+ \sum_{m=1}^{n} \left ( \begin{array}{c} n \\ m \end{array} \right )
\end{array}
\]
which explains that a combination of $n$ loads leads to the sum  of $n$ binary combinations where in each combination $m=1,2,\ldots,n$ loads can be turned ``on'' in
\[
\left ( \begin{array}{c} n \\ m \end{array} \right ) = \frac{n!}{m!(n-m)!}
\]
different ways. Without a required minimum number of $N_i^{\mbox{\tiny on}}$ for the loads, the same number of $2^n$ remains available leading to $2^n\cdot2^n$ combinations at the next time step $t_{k+1}$. 

With $N_i^{\mbox{\tiny on}}>1$, for each combination where $m$ loads have turned ``on'', only $2^{n-m}$ combinations of ``off'' loads are available. This leads to
\[
\sum_{m=0}^n 2^{n-m} \left ( \begin{array}{c} n \\ m \end{array} \right )
\]
at the next time step $t_{k+1}$ for $N_i^{\mbox{\tiny on}}>1$. With
\[
\begin{array}{rcl}
2^n\cdot2^n &=& \displaystyle  \sum_{m=0}^{n} \left ( \begin{array}{c} n \\ m \end{array} \right ) \cdot \sum_{m=0}^{n} \left ( \begin{array}{c} n \\ m \end{array} \right ) \\ &>& \displaystyle  \sum_{m=0}^n 2^{n-m} \left ( \begin{array}{c} n \\ m \end{array} \right )
\end{array}
\]
it is clear that the number of combinations is much smaller. Continuing this argument for the subsequent time steps until $k=N$ emphasizes that number of binary elements in the set ${\cal W}$ is always (much) smaller than $(2^n)^{N-1}$. This results shows that constraints on the allowable load switching helps to alleviate the combinatorial problem, making an  optimization with binary switching computationally feasible.

As an example, consider the case of $n=3$ loads over a power prediction horizon of $N=6$ samples. Without any requirements on minimum number of on/off samples one would have to evaluate $(2^n)^{N-1} = 32768$ possible combinations of the load switching signal $w_i(t_k) \in \{0,1\}$. Starting at a binary combination with all loads off, e.g., $w(0)=[0~0~0]$ and requiring the loads to stay on/off for at least 4 samples reduces the number of possible binary combinations to only $2197$. Clearly, the number of combinations reduces even further for a non-zero initial condition, e.g., $w(0)=[1~0~0]$, where the first load that is switched on is required to stay on over the prediction horizon. An example of this situation is also analyzed in the following figures that show the admissible switching combinations in the graph representation of a tree.

\begin{figure}[ht]\centering
\includegraphics[width=.8\columnwidth]{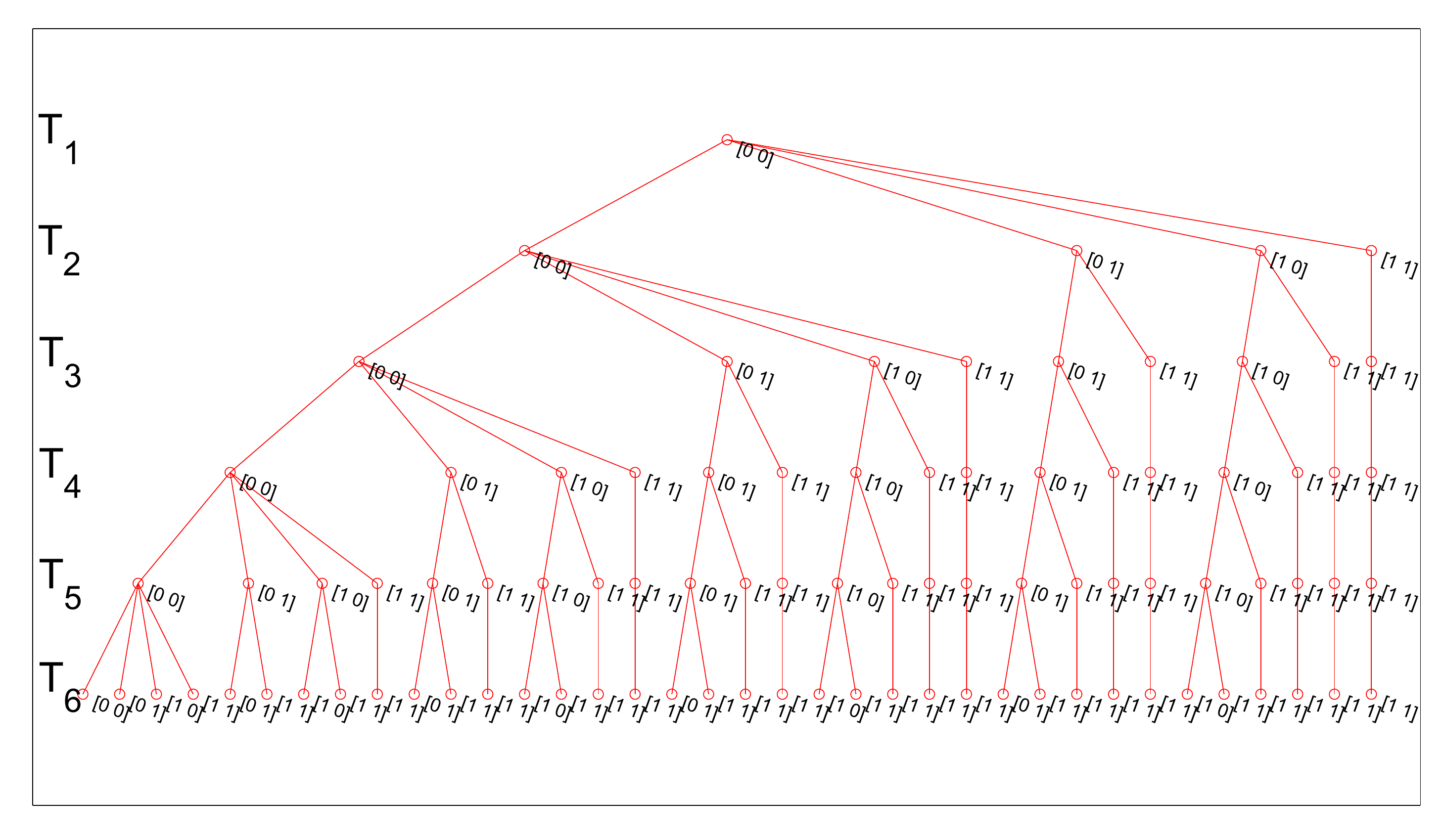}
\caption{Example of tree plot for of $n=2$ loads over a power prediction horizon of $N=6$ with initial conditions $W(\tau_1) = [0~0]$.}
\label{tree1}
\end{figure}

\begin{figure}[ht]\centering
\includegraphics[width=0.8\columnwidth]{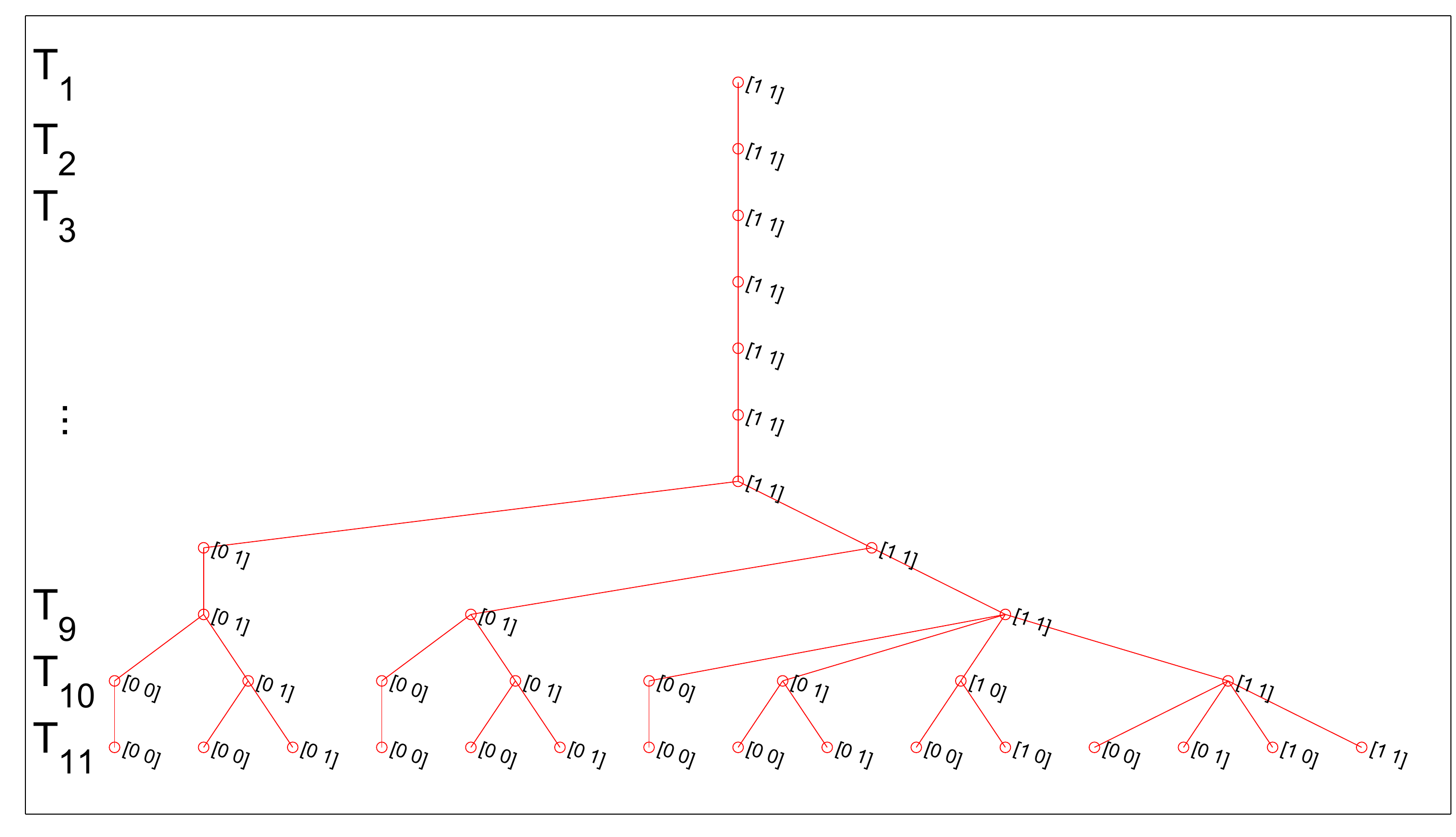}
\caption{Example of tree plot for of $n=2$ loads over a power prediction horizon $N=11$ with initial conditions $W(\tau_1) = [1~1]$.}
\label{tree2}
\end{figure}

In figure \ref{tree1}, a tree plot was created to illustrate the idea, the initial condition here was $[0~0]$ of $n=2$ loads over a power prediction horizon of $N=6$. The total number of combination will be 36. Where the minimum time on or off is 7 and 5 for unit 1 and 2. Another illustrative example is when the initial conditions are $[1~1]$ for a prediction horizon $N=11$ which is large (or long) enough to allow units to be off again, as shown in figure \ref{tree2}.

It is clear from the above illustrations that the number of admissible binary combinations of $n$ loads over a prediction horizon of $N$ points is in general much smaller than $(2^n)^{N-1}$, making the optimization with binary switching computationally feasible for real-time operation. Only a much smaller number of function evaluations have to be performed to find the optimal switching combination over a given prediction horizon, as explained in the following.

\subsection{Moving Horizon Formulation}

Following the power tracking error defined in \fref{powererror}, the dynamic load scheduling optimization problem is formulated as a moving horizon optimization problem
\begin{equation}
\begin{array}{c}w_i(t_m) \\ i=1,2,\ldots,n \\ 
m=k,\ldots,k+N-1\end{array} = \mbox{arg} \min_{w_i(t_m) \in {\cal W}}  f(e(t_l)),
\\
\\
\label{eq:MPC}
\end{equation}
where $l=k,\ldots,k+N-1$ over the admissible set ${\cal W}$, defined in \fref{Wset}.  Similar to the ideas in Model Predictive Control (MPC), the $N \times n$ dimensional optimal switching signal $w_i(t_m)$ is computed over the optimization horizon $m=k,\ldots,k+N-1$. Once the optimal switching signal $w_i(t_m) \in {\cal W}$, $m=k,\ldots,k+N-1$ is computed, the optimal signal is applied to the loads {\em only\/} at the time instant $t_k$, after which the time index $k$ is incremented and the optimization in \fref{MPC} is recomputed over the moving time horizon. 

It should be noted that the admissible set ${\cal W}$ defined in \fref{Wset} has a finite and countable number of binary combinations for the switching signal. Therefore, the $N \times n$ dimensional optimal switching signal $w_i(t_m) \in {\cal W}$ is computed simply by a finite number of evaluation of the criterion function $f(e(t_l))>0$. Hence, no (gradient) based or Quadratic Programming (QP) optimization is used to compute the final value for $w_i(t_k)$ and this makes the problem computationally feasible, even for varying switching load dynamics. Possible candidate functions $f(e(t_l))>0$ may include a least squares criterion
\[
f(e(t_l)) =\sum_{l=k+1}^{k+N} tr\{e(t_l) e(t_l)^T\}
\]
or may include a barrier function 
\begin{equation}
f(e(t_l)) =\sum_{l=k+1}^{k+N} tr\{e(t_l) e(t_l)^T\} - \ln(c(e(t_l)))
\label{eq:optcrit}
\end{equation}
to enforce a positive constraint $e(t_l)>0$. Such constraints may be required to guarantee that the load demand is always smaller than the (predicted) power profile $P(t_k)$ in \fref{powererror}. In this paper we use the quadratic function with a barrier function in \fref{optcrit} to perform tracking of predicted solar power curves by dynamic load switching.

\section{Application Example}\label{sec:AEx}
 
\subsection{Solar Power Data and Load Dynamics}

This model can be implemented for any kind of standalone system (wind, solar or even hybrid) with a forecast tool providing input data. Here we present a standalone solar system connected with 3 units of normalized sizes rated as $x_i=$  60\%, 25.86\% and 12.22\% of full power at a clear day solar curve. San Diego solar power data was collected from UC San Diego campus for illustration purposes of the load scheduling, as solar data variability is prevalent during the summer months when the marine layer clouds are present and leading to multiple days with overcast conditions in the morning and partly cloudy conditions in the afternoon. 
\begin{figure}[ht]\centering
\includegraphics[width=.7\columnwidth]{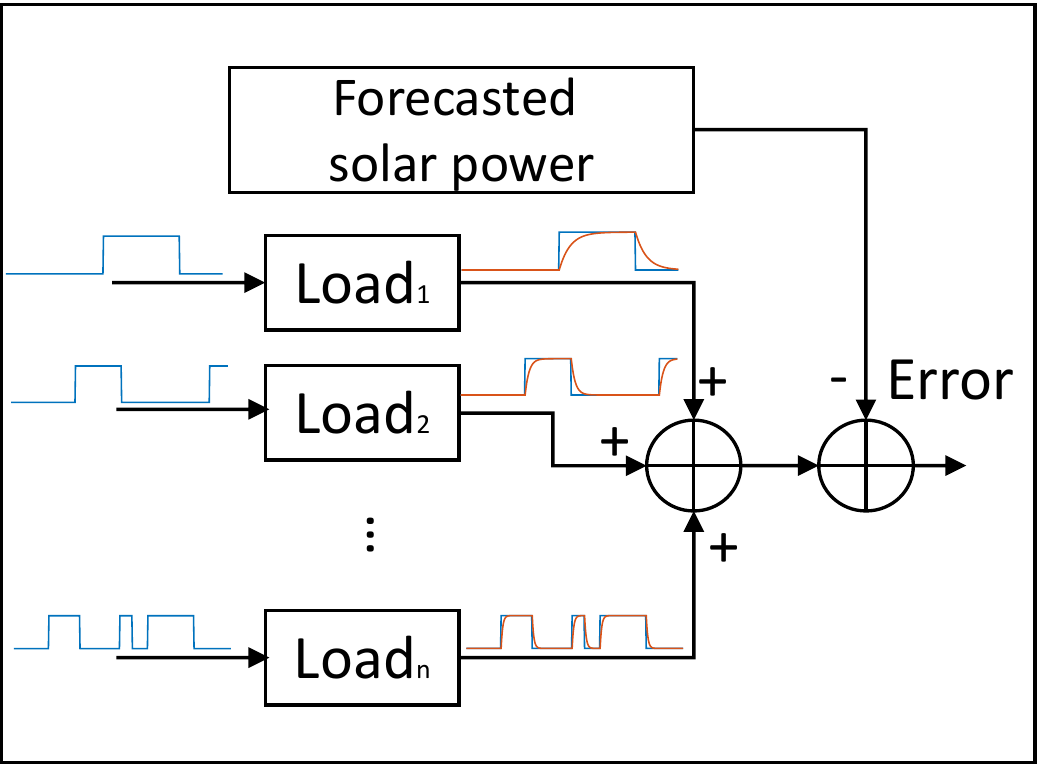}
\caption{Scheduling approach for on/off time of loads based on forecasted solar power.}
\label{systemarch}
\end{figure}
As illustrated in Figure~\ref{systemarch}, the proposed scheduling approach targets to schedule the on/off time of the three different loads to minimize the difference between predicted power delivery and power consumed by the loads, so as to decrease the energy losses as much as possible. For scheduling purposes we also consider every load to have different dynamics for on/off switching modeled by either a first or second system $G_{i}^{on}(s)$ and $G_{i}^{off}(s)$ given earlier in \fref{odeon} and \fref{odeoff}. The dynamics of $G_{i}^{on}(s)$ and $G_{i}^{off}(s)$ is adjusted based on the actual size $x_i$ of the load, assuming that larger loads will take longer to settle to their steady state power demand. A summary of the chosen load dynamic models is summarized in Table~\ref{table1}.
\begin{table}[h]
\centering
\caption{Loads characteristics (seconds)}
\label{table1}
\begin{tabular}{|c|c|c|c|c|c|}
\hline
\begin{tabular}[c]{@{}c@{}}Loads \\ char.\end{tabular} & Size (\%) & Poles$^{on}_i$             & Poles$^{off}_i$ & $T_i^{on}$ & $T_i^{off}$ \\ \hline
$x_1$                                                         & 60.00        & -0.01                & -0.04     & 180        & 180         \\ \hline
$x_2$                                                         & 25.86       & - 0.05 $\pm$ j0.06 & -0.05     & 240        & 240         \\ \hline
$x_3$                                                         & 12.22        & -0.02              & -0.02   & 300        & 300         \\ \hline
\end{tabular}
\end{table}
\subsection{Dynamic Load Switching Results}
\begin{figure}[ht]\centering
\includegraphics[width=.85\columnwidth]{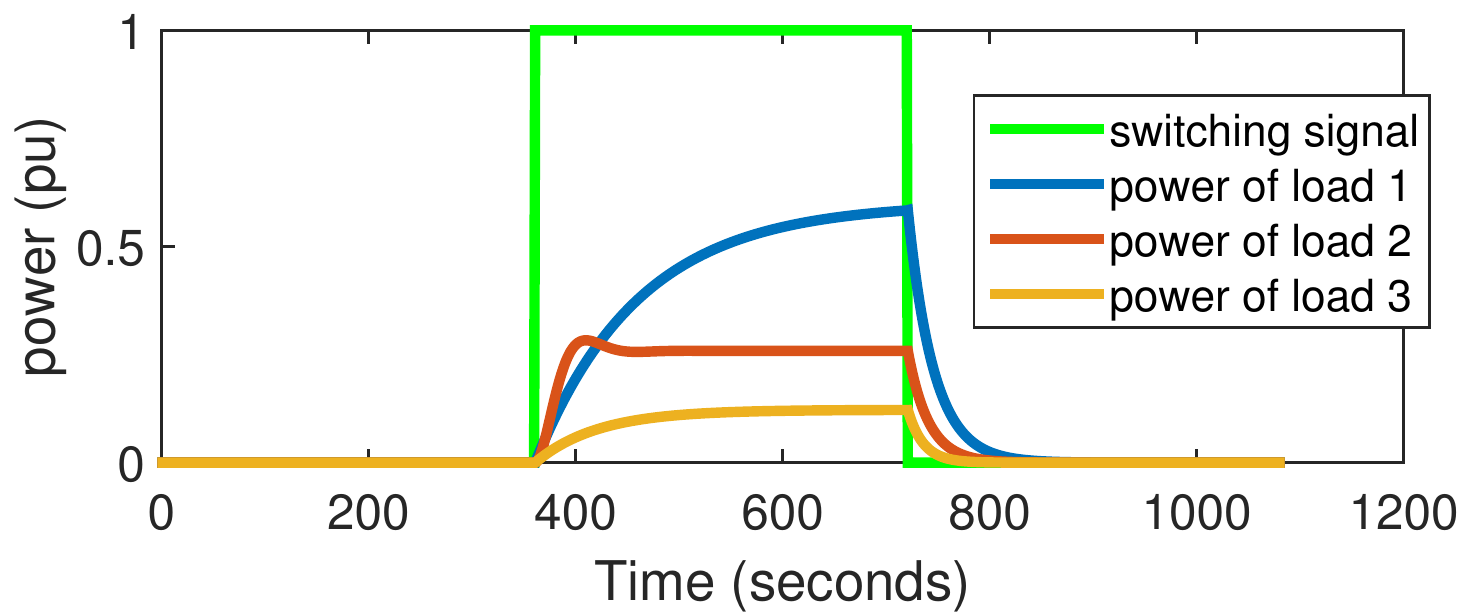}
\caption{Dynamics of the three loads used in this case study, where the dynamics of the loads were summarized in Table~\ref{table1}.}
\label{dynamicsload}
\end{figure}

To illustrate the variability in the dynamics of the loads summarized in Table~\ref{table1}, first the dynamic response of the three loads in our case study are depicted in Figure~\ref{dynamicsload}. it can be seen that Load 1 and Load 3 exhibit first order dynamic behavior, but the dynamics of Load 3 for turning on is much slower due to the larger size of Load 3. Load 2 shows the typical behavior of a second order dynamics model with an initial larger peak load demand for turning on Load 2, e.g., HVAC systems and refrigerators. All loads exhibit similar dynamics when turning off.

\begin{figure}[ht]\centering
\includegraphics[width=.8\columnwidth]{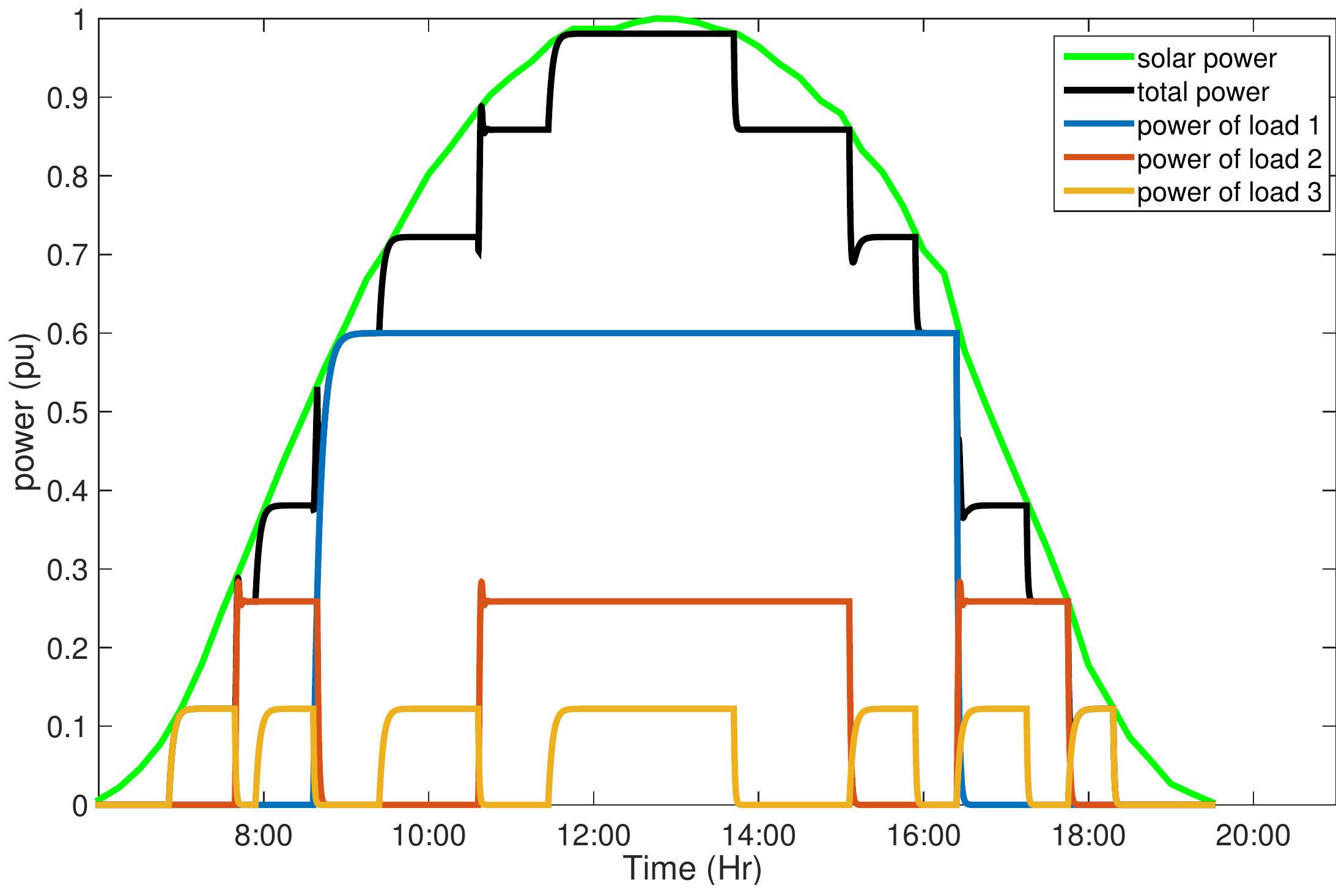}
\caption{Example of load scheduling for $n=3$ loads with size dimensions $x_i=60\%, 25.86\%$ and $12.22\%$ of full solar power, subjected to solar forecast of a stand alone PV on a clear day with power output normalized to peak at 1.}
\label{clearday}
\end{figure}

Using the dynamics of the loads given in Table~\ref{table1} and illustrated in Figure~\ref{dynamicsload} with a sampling time $\Delta_t=1$~second  and a switching time $\tau_i = N_i \Delta_t = 60$~seconds, the three loads are scheduled optimally using a moving horizon of $N=360$~seconds with 6 possible switching times. The results in Figure~\ref{clearday} shows a smooth scheduling for $n=3$ loads over a clear day. Here the solar power is normalized to 1 as well as the loads. It can be seen that the scheduling optimization emphasizes that turning on the largest unit is the main priority. Optimal power tracking is obtained by following the remaining power which is filled in by the remaining smaller loads. During the peak of the day power was too small to keep the smallest unit on.

\begin{figure}[ht]\centering
\includegraphics[width=.8\columnwidth]{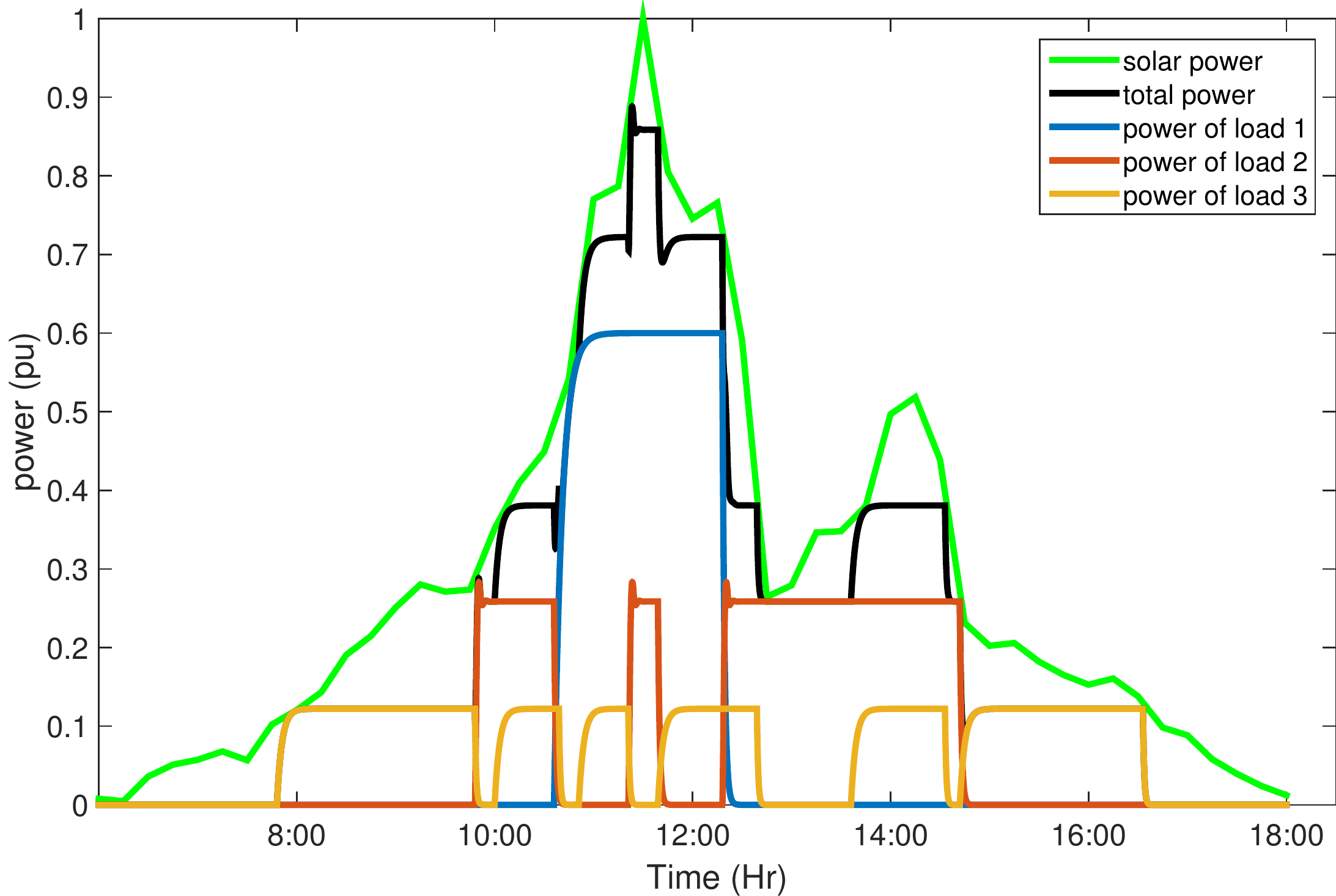}
\caption{Example of load scheduling for $n=3$ loads with size dimensions $x_i=60\%, 25.86\%$ and $12.22\%$ of full solar power, subjected to a solar prediction from a highly variable (partly cloudy) day.}
\label{loadswitching}
\end{figure}

Scheduling for solar power with cloud-induced variability  is a more interesting problem. The results in Figure~\ref{loadswitching} show how (dynamic) scheduling can be done in case of high solar power variability. In this example the smallest load dynamic is assumed to have very short time constants to allow fast power up/down to facilitate fast power tracking. The larger the load the slower the dynamics. Similar scheduling results can also be obtained in the case where loads are assumed to have much larger time constants. Moreover, the computational time was around 45 seconds performed in a 3.4 GHz Intel Core i7 processor with 32 GB of RAM for a full solar day data sampled at 15 minutes. 
\section{Conclusions}\label{sec:CCLS}
Dynamic load scheduling is defined in this paper as the optimal on/off combinations and timing of a set of distinct electric loads via the computation of an optimal binary control signal. For binary load switching, the optimization problem for dynamic load scheduling becomes untractable due to a combinatorial problem where the number of binary combinations grows exponentially with the length of the prediction horizon and the number of loads. In this paper it is shown that constraints on the allowable load switching (typically in the form of a minimum on/off time) help to alleviate the combinatorial problem, making an optimization with binary switching computationally feasible. The model	 presented in this paper is able to solve the on/off scheduling of electric loads with non-linear dynamics based on predictions of the power delivery of a stand-alone solar power source. A prediction of solar power output is used as an input to optimally select the timing and the combinations of a set of given electric loads. The approach is illustrated on electric loads with different dynamics for on versus off switching and solar forecasting data obtained from the Solar Resource Assessment \& Forecasting Laboratory at UC San Diego.
\bibliographystyle{IEEEtran}
\bibliography{mylibACC.bib}
\end{document}